\documentclass[12pt]{article}

\usepackage[table]{xcolor}
\usepackage{amsmath, amsthm}
\usepackage{amsfonts,amssymb}
\usepackage{hyperref}
\usepackage{authblk}
\usepackage[all]{xy}
\usepackage{enumitem}
\usepackage{stmaryrd}
\usepackage[makeroom]{cancel}
\usepackage{bbm}
\usepackage{fancyvrb}
\usepackage[linesnumbered,lined,boxed]{algorithm2e}
\usepackage{xstring}
\usepackage{tikz}
\usetikzlibrary{intersections}

\newcommand{\N}{\mathbb{N}}
\newcommand{\Z}{\mathbb{Z}}
\newcommand{\Q}{\mathbb{Q}}
\newcommand{\R}{\mathbb{R}}

\newcommand{\F}{\mathbb{F}}
\newcommand{\E}{\mathcal{E}}

\renewcommand{\a}{\mathfrak{a}}

\newcommand{\m}{\mathfrak{m}}
\newcommand{\p}{\mathfrak{p}}
\renewcommand{\P}{\mathbb{P}}

\renewcommand{\O}{\mathcal{O}}
\newcommand{\CC}{\mathcal{C}}
\newcommand{\JJ}{\mathcal{J}}
\newcommand{\Gal}{\operatorname{Gal}}

\newcommand{\GL}{\operatorname{GL}}
\newcommand{\SL}{\operatorname{SL}}
\newcommand{\Sp}{\operatorname{Sp}}
\newcommand{\GSp}{\operatorname{GSp}}

\newcommand{\Fl}{{\F_\ell}}

\newcommand{\Ker}{\operatorname{Ker}}

\newcommand{\disc}{\operatorname{disc}}

\newcommand{\Frob}{\operatorname{Frob}}

\renewcommand{\H}{\operatorname{H}_\text{\'et}}

\renewcommand{\S}{\mathcal{S}}

\numberwithin{equation}{section}

\newtheorem{thm}[equation]{Theorem}

\theoremstyle{definition}
\newtheorem{de}[equation]{Definition}
\newtheorem{rk}[equation]{Remark}

\usepackage{float}
\newfloat{Figure}{H}{Fig}[section]
\newfloat{Algorithm}{H}{Algo}[section]
\newfloat{Strategy}{H}{Strat}[section]
\newfloat{Table}{H}{Tab}[section]
\usepackage{fancybox,framed}

\makeatletter
\newcommand{\subjclass}[2][2020]{%
  \let\@oldtitle\@title%
  \gdef\@title{\@oldtitle\footnotetext{#1 \emph{Mathematics subject classification:} #2}}%
}
\newcommand{\keywords}[1]{%
  \let\@@oldtitle\@title%
  \gdef\@title{\@@oldtitle\footnotetext{\emph{Key words and phrases.} #1.}}%
}
\let\c@table\c@equation
\let\c@figure\c@equation
\makeatother
\makeatother

\title{Explicit computation of Galois representations \\ occurring in families of curves}
\author{Nicolas Mascot\thanks{\href{mailto:mascotn@tcd.ie}{mascotn@tcd.ie}}}
\affil{\scriptsize{Trinity College Dublin}}

\VerbatimFootnotes

\begin{document}

\maketitle

\begin{abstract}
We extend our method to compute division polynomials of Jacobians of curves over~$\Q$ to curves over~$\Q(t)$, in view of computing mod~$\ell$ Galois representations occurring in the \'etale cohomology of surfaces over~$\Q$. Although the division polynomials which we obtain are unfortunately too complicated to achieve this last goal, we still obtain explicit families of Galois representations over~$\P^1_\Q$, and we study their degeneration at places of bad reduction of the corresponding curve.
\end{abstract}

\renewcommand{\abstractname}{Acknowledgements}
\begin{abstract}
The author thanks Jean Gillibert for setting him on track to understanding the material presented in Section~\ref{sect:fibres}.
Experiments presented in this paper were carried out using the~\cite{Plafrim} experimental testbed, supported by Inria, CNRS (LABRI and IMB), Universit\'e de Bordeaux, Bordeaux INP, and Conseil R\'egional d’Aquitaine (see \url{https://www.plafrim.fr/}), and on the Warwick mathematics institute computer cluster provided by the EPSRC Programme Grant EP/K034383/1 ``LMF: L-Functions and Modular Forms''. The computer algebra packages used were~\cite{gp} and~\cite{Magma}.
\end{abstract}


\newpage

\section{Introduction}\label{sect:intro}

Suppose we are given a surface~$S$ defined over~$\Q$ as well as a prime~$\ell \in \N$ such that the \'etale cohomology space~$\H^2(S_{\overline \Q},\Z/\ell\Z)$ contains a Galois-submodule which affords a mod~$\ell$ Galois representation~$\rho$ that we wish to compute explicitly. By this, we mean computing a polynomial which encodes~$\rho$ in the following sense:

\begin{de}\label{de:pol_encode_rep}
Let~$K$ be a number field, and let~$\rho : \Gal(\overline K/K) \longrightarrow \GL(V_\rho)$ be a mod~$\ell$ Galois representation, where~$V_\rho$ is an~$\Fl$-vector space of finite dimension. We say that a separable polynomial~$F(x) \in K[x]$ \emph{encodes~$\rho$} if we are given an explicit bijection between~$V_\rho \setminus \{ 0 \}$ and the roots of~$F(x)$ in some extension~$\Omega$ of~$K$ over which~$F(x)$ splits completely, in such a way that the Galois action on the roots of~$F(x)$ matches that on~$V_\rho$. In particular, the splitting field of~$F(x)$ then agrees with the number field~$\overline K^{\Ker \rho}$ cut out by~$\rho$.
\end{de}

In~\cite[2]{SL3}, we sketched a method to compute~$\rho \subset \H^2(S_{\overline \Q},\Z/\ell\Z)$ based on \emph{d\'evissage}~\cite[3.4]{Deligne}, and which may be informally summarised as follows. Pick a proper dominant morphism~$\pi : S \longrightarrow B$ from~$S$ to a curve~$B$ over~$\Q$, and write~$S_b$ for the fibre of~$\pi$ at a point~$b \in B$. Roughly speaking, the Leray spectral sequence~\cite[12.7]{Milne}
attached to~$\pi$ then shows that~$\H^{2} (S_{\overline \Q}, \Z/\ell \Z)$ is made up of~$\H^p\big(B_{\overline \Q}, \H^q(S_b, \Z/\ell \Z)\big)$ for~$p+q=2$. Since the terms for~$p=0,q=2$ and for~$p=2,q=0$ consist of uninteresting bits, we can expect that~$\rho$ occurs in~$\H^1\big(B_{\overline \Q}, \H^1(S_b, \Z/\ell \Z)\big)$. As~$B$ and the~$S_b$ are curves, and as the~$\H^1$ of a curve is essentially the torsion of its Jacobian (see the first part of Theorem~\ref{thm:devissage} below for a precise statement), it is thus reasonable to hope to compute compute~$\rho \subset \H^{2} ({S}_{\overline \Q}, \Z/\ell \Z)$ by:

\begin{Strategy}
\begin{framed}
\begin{enumerate}
\item Computing the family of Galois representations parametrised by~$b \in B$ afforded by the~$\ell$-torsion of the Jacobian of the fibre~$S_b$,
\item Gluing these data into an explicit model of a cover~$C \longrightarrow B$ of curves,
\item Catching~$\rho$ in the~$\ell$-torsion of the Jacobian of the curve~$C$.
\end{enumerate}
\end{framed}
\caption{Computing in the~$\H^2$ of surfaces by looking at the torsion of Jacobians of curves.}
\label{strat:devissage}
\end{Strategy}

The situation is illustrated on Figure~\ref{fig:devissage}.

\begin{Figure}[H]
\begin{center}
\begin{tikzpicture}[scale=2.5]
\draw[thick,variable=\t,domain=-1:1,samples=50]
  plot ({1.9*\t-0.2},{(\t*\t*\t-\t)/5-1}) 
  plot ({2.1*\t+0.2},{(\t*\t*\t-\t)/5+1}) 
  plot ({\t*\t*\t/10-2},{\t}) 
  plot ({(4*\t*\t*\t-\t)/10+2},{\t}); 
\draw(2.2,0) node{$S$};

\draw[thick,variable=\t,domain=-1:1,samples=50]
  plot ({2*\t},{(\t*\t*\t-\t)/5-1.8}); 
\draw(2.2,-1.8) node{$B$};
\draw[thick,->] (0,-1.2) -- (0,-1.55); 
\draw (0.2,-1.375) node{$\pi$};

\draw[gray,variable=\t,domain=-1:1,samples=50] 
  plot ({(1*\t*\t*\t-\t)/10-1},{\t+0.08})
  plot ({(2*\t*\t*\t-\t)/10},{\t})
  plot ({(3*\t*\t*\t-\t)/10+1},{\t-0.08});
  
\draw[gray] (-1.2,1.2) -- (-0.8,1.5) -- (-0.8,3.2) -- (-1.2,2.9) -- (-1.2,1.2);
\draw[gray] (-0.2,1.2) -- (0.2,1.5) -- (0.2,3.2) -- (-0.2,2.9) -- (-0.2,1.2);
\draw[gray] (0.8,1.2) -- (1.2,1.5) -- (1.2,3.2) -- (0.8,2.9) -- (0.8,1.2);

\fill[red] (-1.05,2.19) circle (1pt);
\fill[red] (-1.02,2.9) circle (1pt);
\fill[red] (-0.97,2.2-0.33) circle (1pt);
\fill[red] (0,2.2) circle (1pt);
\fill[red] (0.05,2.705) circle (1pt);
\fill[red] (0.025,1.68) circle (1pt);
\fill[red] (0.99,2.2) circle (1pt);
\fill[red] (1.01,2.53) circle (1pt);
\fill[red] (0.95,1.51) circle (1pt);

\draw[thick,red,variable=\t,domain=-1.2:1.2,samples=100]
  plot ({(1.5+\t*\t/3)*(sin(200*\t)+sin(600*\t)/20)},{2.2+1.8*\t-atan(1.6*\t)/50+cos(500*\t)/400});
\draw(1.65,2.32) node{$C$};

\end{tikzpicture}
\end{center}
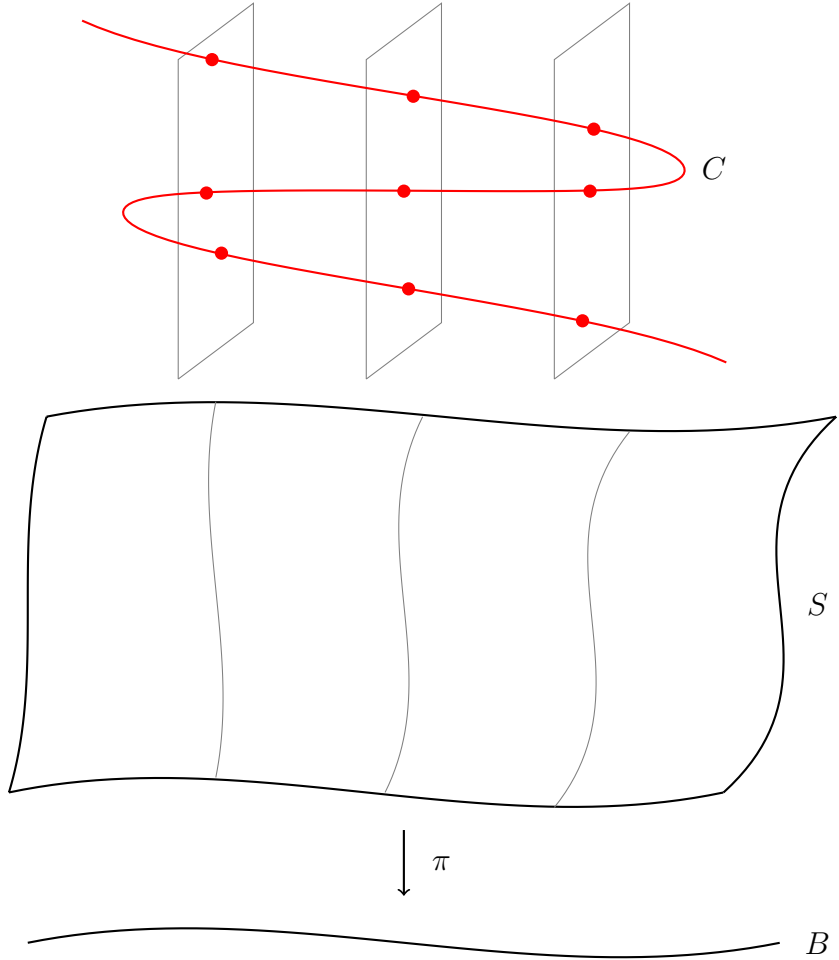
\caption{The surface~$S$ with some of the fibres~$S_b$ of~$\pi$. The rectangles above them represent the Jacobian of these fibres, inside which the red dots represent~$\ell$-torsion points. These points define a curve~$C$ whose Jacobian should contain~$\rho$ in its~$\ell$-torsion.}
\label{fig:devissage}
\end{Figure}

More precisely, we have the following result:

\begin{thm}\label{thm:devissage}
Given an~$\Fl$-Galois-module~$M$ and an integer~$n \in \Z$, write~$M(n)$ for the twist of~$M$ by the~$n$-th power of the mod~$\ell$ cyclotomic character.
\begin{enumerate}
\item Let~$X$ be a nonsingular, geometrically irreducible curve over  a number field~$K$, and let~$J$ be the Jacobian of the completion of~$X$. If~$X$ is complete, then~$\H^1(X_{\overline K},\Z/\ell \Z) \simeq J[\ell](-1)$ as Galois modules. If~$X$ is not complete, then~$\H^1(X_{\overline K},\Z/\ell \Z)$ is an extension of~$J[\ell](-1)$ by copies of~$(\Z/\ell\Z)(-1)$.
\item Suppose~$\rho$ is a mod~$\ell$ Galois representation contained in~$\H^2(S_{\overline \Q},\Z/\ell\Z)$ (up to semi-simplification). Let~$B' = B \setminus Z$, where~$Z \subset B$ is the locus of bad fibres of~$\pi$. Assume that~$\rho$ has no Jordan-Hölder components of the form~$(\Z/\ell\Z)(n)$ for any~$n \in \Z$, and no component in common with~$\eta(-1)$, where~$\eta$ is the mod~$\ell$ permutation representation induced by the Galois action on the geometrically irreducible components of the bad fibres of~$\pi$. Then~$\rho$ is also contained (up to semi-simplification) in~$\H^1(C_{\overline \Q},\Z/\ell \Z)(-1)$, where~$C$ is the completion of the cover of~$B'$ formed by the nonzero~$\ell$-torsion points of the Jacobian of the~$S_b$.
\end{enumerate}
\end{thm}

Part 1 is standard (cf.~\cite[14.2,14.4,16.2]{Milne}), and part 2 is~\cite[Thm 7]{SL3}. In particular, if~$\rho$ satisfies the assumptions of part 2, and if~$C$ is geometrically irreducible, then~$\rho$ is found (up to twist) in the~$\ell$-torsion of the Jacobian of~$C$. More generally, if~$C$ is not geometrically irreducible, consider a Galois number field~$K \subset \overline \Q$ such that the geometrically irreducible components~$C_i$ of~$C$ are defined over~$K$; then~$\rho$ will be found in the induction to~$\Gal(\overline \Q / \Q)$ of the representation of~$\Gal(\overline \Q/K)$ afforded by the~$\ell$-torsion of the Jacobians of the~$C_i$.

Let us now explain in more detail how to turn these observations into an algorithm to compute~$\rho$ explicitly, assuming for simplicity that~$C$ is geometrically irreducible. In~\cite{Hensel}, we described an algorithm which, given a proper, nonsingular, and geometrically irreducible curve~$C$ over a number field\footnote{At present, this algorithm is only implemented for~$K=\Q$, but its generalisation to number fields is straightforward.}~$K$ and a prime~$\ell \in \N$, computes what may be called an~$\ell$-division polynomial~$R_{C,\ell}(x) \in K[x]$ of~$C$, that is to say a polynomial which encodes the representation afforded by the~$\ell$-torsion of the Jacobian~$J$ of~$C$ in the sense of Definition~\ref{de:pol_encode_rep}. This algorithm is also capable of computing the subrepresentation afforded by a Galois-submodule~$V$ of~$J[\ell]$, provided that there exists a prime~$\p \nmid \ell$ of~$K$ where~$C$ has good reduction and such that~$V \subset J[\ell]$ may be characterised by the characteristic polynomial of~$\Frob_\p$ acting on~$V$.

Suppose for the sake of the exposition that we are given an equation~$f(x,y,t) \in \Q[x,y,t]$ such that our surface~$S$ is the desingularisation of the projective closure of the patch defined by~$f(x,y,t)=0$. It is then natural to choose~$B = \P^1_\Q$ and~$\pi$ the projection~$(x,y,t) \mapsto t$, thereby viewing the surface~$S$ as a curve~$\S$ over~$\Q(t)$. Suppose furthermore we generalised our division polynomial algorithm~\cite{Hensel} to curves over~$\Q(t)$. We would then be able to compute a division polynomial~$R_{\S,\ell}(x,t) \in \Q(t)[x]$ for~$\S$, whose specialisation~$R_{\S,\ell}(x,t_0) \in \Q(t_0)[x]$ at any good fibre~$t=t_0 \in B$ of~$\pi$ would an~$\ell$-division polynomial of the fibre~$S_{t_0}$. Then the equation~$R_{\S,\ell}(x,t)=0$ would define the curve~$C$ such that~$\rho$ occurs (up to twist by the cyclotomic character) in the~$\ell$-torsion of the Jacobian of~$C$, so that we may compute~$\rho$ by applying the original version of~\cite{Hensel} to~$C$, by isolating the twist of~$\rho$ in the Jacobian~$J_C$ of~$C$ from the knowledge of the characteristic polynomial of~$\rho(\Frob_\p)$ where~$\p$ is as described above (cf.~\cite{SL3} for a successfully worked out example of this approach).

In particular, we would not even need to compute all of the~$\ell^{2g_C}$ points of~$J_C[\ell]$, which would be impractical even for~$\ell=2$ as soon as the genus~$g_C$ of~$C$ is moderately large, but only the~$\ell^{\deg \rho}$ points of the subspace affording the twist of~$\rho$ contained in~$J_C[\ell]$. On the other hand, this method forces us to compute all the~$\ell$-torsion points of the Jacobian of~$\S$ in order to get an equation for~$C$, and this therefore only applicable when the genus of~$\S$ is reasonably small.

The purpose of this article is to explain how~\cite{Hensel} can indeed be generalised to curves over~$\Q(t)$, thereby making it theoretically possible to compute explicitly mod~$\ell$ Galois representations which occur in the~$\H^2$ of surfaces.

\begin{rk}
Very general but unfortunately impractical algorithms to compute with \'etale cohomology are presented in~\cite{Madore} and~\cite{Testa}. In contrast, our goal is to obtain a practical method for the specific case of the~$\H^2$ of surfaces.
\end{rk}

We show how~\cite{Hensel} can be generalised to curves over~$\Q(t)$ in Section~\ref{sect:HenselQt}. Since~\cite{Hensel} requires the curve to be given as a Riemann-Roch space, in Section~\ref{sect:algcurves} we briefly recall how to perform various computations with plane algebraic curves, including the determination of Riemann-Roch spaces and the verification whether the curve is geometrically irreducible.

As an application, in Section~\ref{sect:examples} we compute division polynomials~$R_{\S,\ell}(x,t)$ for three curves~$\S$ over~$\Q(t)$, of respective genera~$1$,~$2$, and~$3$. This makes it possible, in principle, to compute with the~$\H^2$ of the corresponding surfaces over~$\Q$; but unfortunately, the equations which we obtain for the curves of genera~$2$ and~$3$ are too complicated for this to practical. However, the data that we obtain is still worth our attention, since it encodes families of Galois representations over~$B = \P^1_\Q$, and it is especially interesting to study how these families degenerate at bad fibres, which we do in Section \ref{sect:degen}; in particular, we strive to find a geometric explanation for the ramification of these degenerations.
%


\section{Division polynomials over~$\Q(t)$}\label{sect:HenselQt}

\subsection{Sketch of the algorithm over~$\Q$}

Let still~$\ell \in \N$ be prime. The purpose of this section is to explain how our algorithm~\cite{Hensel} to compute~$\ell$-division polynomials of curves over~$\Q$ can be generalised to curves over~$\Q(t)$. In this view, let us first recall how this algorithm works with a curve~$C$ over~$\Q$:

\begin{Algorithm}
\begin{framed}
\begin{enumerate}
\item Pick a prime~$p \neq \ell$ of good reduction of~$C$. Determine~$a \in \N$ such that the~$\ell$-torsion of the Jacobian~$J$ of~$C$ is defined over~$\F_{q}$, where~$q=p^a$.
\item Generate points of~$J(\F_{q})[\ell]$ which span~$J[\ell]$ as an~$\Fl[\Frob_p]$-module.
\item Lift these points to~$J(\Z_q/p^e)[\ell]$, where~$\Z_q$ is the ring of integers of the unramified extension of~$\Q_p$ with residue field~$\F_q$, and~$e \in \N$ is an accuracy parameter.
\item Construct an evaluation map~$\alpha \in \Q(J)$.
\item Expand~$\displaystyle \tilde F(x) = \prod_{0 \neq t \in J[\ell]} \!\!\!\!\! \big(x - \alpha(t)\big) \in (\Z/p^e\Z)[x]$, and identify it as an element~$F(x)$ of~$\Q[x]$.
\end{enumerate}
\end{framed}
\caption{Division polynomial of a curve over~$\Q$.}
\label{alg:divpol_Q}
\end{Algorithm}

The idea is thus to pick an auxiliary prime~$p$, and to rely on the fact that~$J[\ell]$ is \'etale at~$p$ to construct~$p$-adic approximations of points of~$J[\ell]$.

The polynomial~$F(x)$ is then an~$\ell$-division polynomial of~$C$ in the sense of Definition~\ref{de:pol_encode_rep}. This supposes that~$\alpha$ is defined and injective on~$J[\ell]$; if this is not the case, we start over with another~$\alpha$. This also supposes that the accuracy parameter~$e$ is large enough to identify~$F(x)$ from its mod~$p^e$ approximation~$\tilde F(x)$. In particular, the correctness of this method is not rigorously guaranteed, although this could be done by confirming that the elements of~$J(\Z_q/p^e)[\ell]$ are indeed~$p$-adic approximations of~$\ell$-torsion points defined over the stem fields of the irreducible factors of~$F(x)$. Besides, in most cases, one easily convinces oneself beyond reasonable doubt that the output~$F(x)$ is correct, e.g. by checking that it has the appropriate Galois group and ramification.

In order to compute in~$J$, this algorithm relies on Makdisi's algorithms~\cite{Mak1,Mak2}. These algorithms were originally designed to work over a field, so in~\cite{Hensel} we generalised them to work over a local ring such as~$\Z_q/p^e$. These algorithms also require the knowledge of an explicit basis of a Riemann-Roch space of~$C$ of high-enough degree so as to represent~$C$ internally (cf. the bottom of page 1421 in~\cite{Hensel}), so we will explain in Section~\ref{sect:algcurves} below how such a basis may be computed from a (possibly singular) plane model of~$C$.

\subsection{Sketch of the algorithm over~$\Q(t)$}

By analogy with the embedding of~$\Q$ into its completion~$\Q_p$, it is natural to extend Algorithm~\ref{alg:divpol_Q} to curves over~$\Q(t)$ by embedding~$\Q(t)$ into the~$p$-adic Laurent series field~$\Q_p((t))$. This leads to the following idea to compute an~$\ell$-division polynomial of a curve~$\CC$ over~$\Q(t)$:

\begin{Algorithm}
\begin{framed}
\begin{enumerate}
\item If required, shift the parameter~$t$ so that~$\CC$ has good reduction~$C_0$ at~$t=0$. Pick a prime~$p \neq \ell$ of good reduction of~$C_0$, and determine~$a \in \N$ such that the~$\ell$ torsion of the Jacobian~$J_0$ of~$C_0$ is defined dover~$\F_{q}$, where~$q=p^a$.
\item Generate points of~$J_0(\F_{q})[\ell]$ which span~$J_0[\ell]$ as an~$\Fl[\Frob_p]$-module.
\item Lift these points to~$\JJ(R)[\ell]$, where~$\JJ$ is the Jacobian of~$\CC$ and~$R$ is a finite quotient of the formal power series ring~$\Z_q[[t]]$.
\item Construct an evaluation map~$\alpha \in \Q(t)(\JJ)$.
\item Expand~$\displaystyle \tilde F(x) = \prod_{0 \neq t \in J[\ell]}\big(x - \alpha(t)\big) \in R[x]$, and identify it as an element~$F(x)$ of~$\Q(t)[x]$.
\end{enumerate}
\end{framed}
\caption{Division polynomial of a curve over~$\Q(t)$.}
\label{alg:divpol_Qt}
\end{Algorithm}

This assumes that we manage to extend Makdisi's algorithms to finite quotients of~$\Z_q[[t]]$. This is actually not an issue, because the extension which we designed in~\cite{Hensel} works with any finite local ring~$R$ over which one can perform linear algebra in ``good reduction cases'' in the following sense:
\begin{de}\label{def:linalg_goodred}
Let~$R=\O/\a$ be finite quotient of a local domain~$\O$. Let~$K$ be the fraction field of~$\O$, and let~$k$ be the residue field of~$\O$. We say that \emph{we can perform linear algebra over~$R$ in cases of good reduction} if, given the reduction mod~$\a$ of a matrix~$A$ over~$\O$ such that the rank of~$A$ is the same over~$K$ and over~$k$, we can compute an approximation in~$\O/\a$ of a~$K$-basis of the kernel of~$A$.
\end{de}
Similarly, the construction~\cite[2.2.3]{Hensel} of evaluation maps~$\alpha$ generalises to Jacobians of curves over~$\Q(t)$ without change.

Finally, we can identify the coefficients of~$\tilde F(x)$ as elements of~$\Q(t)$ by a combination of~$p$-adic rational reconstruction (as we did in the original version of~\cite{Hensel}) and of Pad\'e approximants (see Remark~\ref{rk:padic_Pade} below for practical details).

\subsection{Lifting torsion points~$(p,t)$-adically}

In order to turn these ideas into a proper algorithm, we still must explain what kind of finite quotients~$R$ of~$\Z_q[[t]]$ we will work with, and how to lift an~$\ell$-torsion point from~$\F_q$ to~$R$.

A first natural choice for~$R$ would be~$R_e = \Z_q[[t]]/\m^e$, where~$\m=(p,t)$ is the maximal ideal of~$\Z_q[[t]]$ and~$e \in \N$ is an accuracy parameter as in Algorithm~\ref{alg:divpol_Q}. This choice may be appealing at first, as it would give us the hope of being able to raise the~$p$-adic and the~$t$-adic accuracy of torsion points simultaneously; but unfortunately, we will see below that~$R_e$ having Krull dimension 2 actually results in an algorithmic obstacle to lifting torsion points. Furthermore, elements of~$R_e$ are of the form~$\sum_{j<e} \lambda_j t^j$ where~$\lambda_j \in \Z/p^{e-j}\Z$ is known with poor accuracy for large~$j$; as a result, in~$\tilde F(x)$, the coefficients of high powers of~$t$ would be known with poor~$p$-adic accuracy, which would force us to increase the value of~$e$ so as to identify them, so we would end up lugging around high powers of~$t$ throughout the calculation only to drop them at the final stage since they are~$p$-adically too imprecise to be identified as rational numbers, and thus result in a major waste of time.

We have therefore decided to work with the quotients~$R=R_{e,h} = (\Z_q/p^e\Z_q)[t]/(t^h)$, where~$h \in \N$ is a second accuracy parameter. The introduction of this new parameter grants us the flexibility of setting the~$p$-adic accuracy independently from the~$t$-adic one, which turns out to be useful in practice. Furthermore, this makes it possible to generalise our algorithm to lift torsion points. In order to see why, recall how we proceeded over~$\Q$ in~\cite{Hensel}:

Let~$\O=\Z_q$,~$\varpi=p$,~$\m = \varpi \O$,~$K=\Q_q$, and let~$J$ be the Jacobian of a curve over~$K$ which has good reduction at~$\m$. Given~$e \in \N$, a point~$x \in J(\O/\m^e)$ is represented in Makdisi's algorithms (as generalised in~\cite{Hensel}) by a matrix~$W_x$ with entries in~$\O/\m^e$; but conversely, most such matrices do not represent any point of~$J$. We thus began with an algorithm~\cite[Algorithm 9]{Hensel} which, given an integer~$e \in \N$ and a matrix~$W_x$ representing~$x \in J(\O/\m^e)$, computes a lift of~$W_x$ to~$\O/\m^{2e}$ which represents a lift of~$x$ to~$J(\O/\m^{2e})$.

Due to the tangent space of~$J$ at~$x$, this lift of~$x$ is not unique, and indeed this algorithm can return several matrices representing different random lifts of~$x$ if required. But this also means that even if~$x$ was~$\ell$-torsion in~$J(\O/\m^e)$, none of these lifts to~$J(\O/\m^{2e})$ are guaranteed (nor even likely) to be~$\ell$-torsion.

In order to circumvent this problem, we showed how to construct an algebraic ``coordinate chart''~$\kappa : U \hookrightarrow \O^n$, where~$n$ is a fixed integer not smaller than the genus~$g$ of the curve. This chart is defined on an~$\m$-adic neighbourhood~$U$ of the origin~$0 \in J(\O)$, and turns the mod~$\m^e$ representation in Makdisi form of a point~$x \in U$ into a vector~$\kappa(x) \in (\O/\m^e)^n$ such that for all~$e' \leqslant e$,~$\kappa(x) = 0 \bmod \m^{e'}$ if and only if~$x=0$ in~$J(\O/\m^{e'})$. As~$\O$ is furthermore principal with uniformiser~$\varpi=p$, we then designed a second algorithm~\cite[Algorithm 11]{Hensel}, which computes the unique lift to~$J(\O/\m^{2e})[\ell]$ of a point~$x \in J(\O/\m^e)[\ell]$ as follows:

\begin{Algorithm}
\begin{framed}
\begin{enumerate}
\item Use algorithm~\cite[Algorithm 9]{Hensel} to generate~$g+1$ matrices~$W_0,\cdots,W_g$ representing random lifts~$x_0,\cdots,x_g$ of~$x$ to~$J(\O/\m^{2e})$.
\item \label{alg:Maklift_2} For each of these lifts, compute the vectors~$\displaystyle k_i = \frac1{\varpi^e} \kappa([\ell]x_i) \in (\O/\m^e)^n$.
\item Try to find scalars~$\lambda_1, \cdots, \lambda_g \in \O/\m^{2e}$ such that~$\sum_{i=0}^g \lambda_i k_i = 0 \bmod \m^e$ and~$\sum_{i=0}^g \lambda_i = 1 \bmod \m^{2e}$, and return the matrix~$\sum_{i=0}^g \lambda_i W_i$.
\end{enumerate}
\end{framed}
\caption{Lifting an~$\ell$-torsion point in Makdisi form.}
\label{alg:Maklift}
\end{Algorithm}

The idea is that with high probability, the lifts~$x_i$ form an affine coordinate frame of the tangent space of~$J$ at~$x$, which guarantees the existence and uniqueness of the~$\lambda_i$ (and otherwise, we start over with other random lifts~$x_i$). Note that since~$x$ is assumed to be~$\ell$-torsion mod~$\m^e$, we have~$\kappa(x_i) = 0 \bmod \m^e$ for all~$i$, so division by~$\varpi^e$ does result in the~$k_i$ being integral. This division is essential so that we can find the~$\lambda_i$ by solving a linear system over the local ring~$\O/\m^{2e}$, since it ensures that this system will have good reduction in the sense of Definition~\ref{def:linalg_goodred} provided as long as the~$x_i$ do form an affine frame.

Let us now see how to generalise Algorithm~\ref{alg:Maklift} to the case where~$\O = \Z_q[[t]$. We can now see why working with quotients of~$\Z_q[[t]]$ of the form~$\Z_q[[t]]/(p,t)^e$ would be an issue: In step~\ref{alg:Maklift_2}, we would obtain vectors~$\kappa([\ell]x_i)$ with entries in~$(p,t)^e/(p,t)^{2e}$, but since the ideal~$(p,t)$ is not principal, we would not be able to renormalise the linear system defining the~$\lambda_i$ into a system of good reduction in the sense of Definition~\ref{def:linalg_goodred}.

In contrast, by working with quotients of the form~$(\Z_q/p^e\Z_q)[t]/(t^h)$, we can generalise Algorithm~\ref{alg:divpol_Q} as follows: given a point~$x \in J_0(\F_q)[\ell]$, we can first lift it~$p$-adically to~$J(\Z_q/p^e)[\ell]$ by using the original version of Algorithm~\ref{alg:Maklift} as described in~\cite{Hensel}, and then, we can lift this lift~$t$-adically to~$\JJ\big((\Z_q/p^e\Z_q)[t]/(t^h)\big)[\ell]$, by applying Algorithm~\ref{alg:Maklift} with~$\O = (\Z_q/p^e\Z_q)[[t]]$ and~$\varpi = t$. Indeed, even though~$\m = t \O$ is no longer maximal, the point is that the quotient~$\O/\m^h = (\Z_q/p^e\Z_q)[t]/(t^h)$ is still a local ring with ``residue ring''~$k = \Z_q/p^e$ which is still local, so that our generalisation of Makdisi to local rings is able to handle working over it.

We are thus able to lift torsion points from~$J_0(\F_q)[\ell]$ to~$\JJ\big((\Z_q/p^e\Z_q)[t]/(t^h)\big)[\ell]$, and thus to extend our method~\cite{Hensel} to curves defined over~$\Q(t)$.

\begin{rk}\label{rk:padic_Pade}
In practice, when we identify elements~$c \in \Q(t)$ from an approximation in~$(\Z/p^e\Z)[t]/(t^h)$ at the end of Algorithm~\ref{alg:divpol_Qt}, rather than first identifying~$c$ as an element of~$\Q[t]/(t^h)$ by~$p$-adic rational reconstruction and then as an element of~$\Q(t)$ by Pad\'e approximants over~$\Q$, it is much more efficient to proceed in the reverse order, that is to say to first use Pad\'e approximants over~$\Q_p$ so as to identify~$c$ as an element of~$\Q_p(t)$ whose coefficients are known mod~$p^e$, and then to reconstruct these coefficients as rational numbers. The reason for this is that unless~$h$ is quite small, the Taylor coefficients of~$c$ up to~$O(t^h)$ will typically have a very large arithmetic height, so that identifying them would require the~$p$-adic precision parameter~$e$ to be very high, which would drastically reduce the execution speed of the whole of Algorithm~\ref{alg:divpol_Qt}. For example, in Section~\ref{sect:divpol_Q} below, identifying the coefficients of a~$2$-division polynomial of a family of plane quartics requires~$h=128$, and experimentation has shown to us that this in turn requires~$e=4096$ with the first method, but only~$e=128$ with the second one. 
\end{rk}

\section{Computing with plane algebraic curves}\label{sect:algcurves}

When we apply Strategy~\ref{strat:devissage}, on both occasions when we use our algorithm to compute an~$\ell$-division polynomial of a curve (first over~$\Q(t)$ with Algorithm~\ref{alg:divpol_Qt}, and then over~$\Q$ with Algorithm~\ref{alg:divpol_Q}), that curve is given to us by a plane equation, which is possibly singular. However, as explained in the previous Section, our~$\ell$-division polynomial algorithm relies on Makdisi's algorithms, which require the curve to be represented by a Riemann-Roch space of high-enough degree.

The purpose of this Section is therefore to explain how one may perform explicit computations, such as Riemman-Roch spaces, with curves given by possibly singular plane models. Such functionalities are already available in some computer algebra packages such as~\cite{Magma}, but our implementation of the~$\ell$-division polynomial algorithm is based on~\cite{gp}, and converting data from~\cite{Magma} to~\cite{gp} is tedious and tends to break the flow of automation. We have therefore implemented our own package to compute with plane algebraic curves in~\cite{gp}, in a way which is tailored towards our needs.

\subsection{Representing the desingularised curve}\label{sect:Duval}

Fix a ground field~$K$ over which one can algorithmically factor polynomials and perform linear algebra. For example,~$K$ could be~$\Q$ or~$\Q(t)$. We also assume that~$K$ has characteristic 0, although this is hypothesis is not essential (see Remark~\ref{rk:Duval_char0} below).

Suppose we are given an irreducible polynomial~$f(x,y) \in K[x,y]$. It defines an affine curve~$C$ over~$K$, but instead one typically wants to work with~$\tilde C$, the desingularisation of the projective completion of~$C$. Nonsingular points of~$C$ may be identified with points of~$\tilde C$, so we only need a specific representation for points of~$\tilde C$ at infinity or above singular points of~$C$.

One possibility would be to construct an explicit model of~$\tilde C$ made up of several charts in a higher-dimensional ambient space; however, this approach would lead to Gr\"obner bases calculations in many variables, which could be very slow. Therefore, we have instead decided to represent these points of~$\tilde C$ by formal series parametrisations. For instance, if~$f(x,y)=xy+\cdots$ so that~$C$ has a node at the origin, the two points of~$\tilde C$ corresponding to the two branches of this node can be represented by parametrisations of the form
\[ x=t, \ y=t+O(t^2) \quad \text{and} \quad x=t, \ y=-t+O(t^2). \]

In order to compute such parametrisations, we can take advantage of the fact that the field~$\overline K\{\{x\}\}$ of Puiseux series over~$\overline K$ contains an algebraic closure of~$K(x)$: for each root~$y = \sum_{m \geqslant m_0} a_m x^{m/e} \in \overline K\{\{x\}\}$ of~$f(x,y) \in K(x)[y]$, we obtain the parametrisation 
\begin{equation} x=t^e, \ y =  \sum_{m \geqslant m_0} a_m t^m \in \overline K((t)). \label{eqn:branch_Puiseux} \end{equation}
One might thus hope for a bijection between the points of~$\tilde C$ above~$x=0$ and parametrisations of the form~$x=t^e$,~$y \in \overline K((t))$ with~$x$ and~$y$ not both series in~$t^m$ for any~$m \geqslant 2$; but unfortunately, this is not the case, because \eqref{eqn:branch_Puiseux} can be reparametrised as
 \[ x=t'^e, \ y =  \sum_{m \geqslant m_0} \zeta^m a_m t'^m \]
 where~$t=\zeta t'$ for any~$e$-th root of unity~$\zeta \in \overline K$. In particular, with this approach, there would be no hope to match the extension of~$K$ generated by the coefficients~$a_j$ with the field of definition of the corresponding point\footnote{Unless of course~$K$ happens to contain the roots of unity of all orders, which typically will not be the case for the applications which we have in mind since we will be working over~$K=\Q$ or~$\Q(t)$.}.
 
Fortunately, Duval~\cite{Duval} has shown that these problems can be circumvented by allowing parametrisations of the form~$x = b t^e$,~$y \in \overline K((t))$ where~$b \in \overline K$ is a constant:
 
 \begin{thm}\label{thm:Duval}
 Let~$f(x,y) \in K[x,y]$ be irreducible of degree~$n$ in~$y$. There exists a finite set of parametrisations
 \[ x=b_j t^{e_j}, \ y = \sum_{m \geqslant m_j} a_{j,m} t^m \]
 where for each~$j$, the~$b_j$ and the~$a_{j,m}$ lie in~$\overline K$ and span a finite extension~$L_j$ of~$K$, and such that the~$n$ roots of~$f$ in~$\overline K \{\{x\}\}$ are obtained without repetition as
 \[ y = \sum_{m \geqslant m_j} a_{j,m}^\sigma (\beta x^{1/e_j})^m \]
 where~$\sigma$ ranges over the~$K$-embeddings of~$L_j$ into~$\overline K$ and~$\beta$ ranges over~$\{ \beta \in \overline K \, \vert \, \beta^{-e_j} = b_j^\sigma \}$ (so that~$t = \beta x^{1/e_j}$ is what one obtains when solving~$x=b_j t^{e_j}$ for~$t$).
 \end{thm}

This means that we have a Galois-equivariant bijection between this set of parametrisations and the set of places of the function field~$K(C) = K(x)[y]/f(x,y)$ of~$C$ above~$x=0$, and therefore with the points of~$\tilde C$ above~$x=0$. In particular, we have
\[ \sum_j e_j f_j = n \]
where the~$f_j = [L_j:K]$ are the residue degrees and the~$e_j$ are the ramification indices, so that the~$L_j$ are the fields of definition of the corresponding points of~$\tilde C$, and that the
\begin{equation} \prod_{\sigma : L_j \hookrightarrow \overline K} \prod_{\beta^{-e_j} = b_j} \left(y -  \sum_{m \geqslant m_j} a_{j,m}^\sigma (\beta x^{1/e_j})^m \right) \end{equation}
are the irreducible factors of~$f(x,y)$ over~$K((x))$. Note the analogy with the determination of the decomposition of a prime number~$p$ in a number field by studying the factorisation over~$\Q_p$ of a polynomial defining that number field.


Duval explains that these parametrisations can be computed as follows:

\begin{Algorithm}
\begin{framed}
\begin{enumerate}
\item\label{alg:Duval_1} Draw the Newton polygon of~$f(x,y)$, that is to say the lower convex hull of the points~$(i,j) \in \Z^2$ such that the coefficient~$a_{i,j}$ of~$y^i x^j$ in~$f(x,y)=\sum_{i,j} a_{i,j} y^i x^j$ is nonzero.
\item For each segment~$pi+qj=r$ of the Newton polygon, where~$p,q,r \in \Z$ and~$\gcd(p,q)=1$, find~$u,v \in \Z$ such that~$up+vq=1$, and let~$f_0 = \sum_{pi+qj=r} a_{i,j} x^j y^i$. Then for each~$b \in \overline K$ such that~$f_0(b^{-u} t^q,b^v t^p)=0$, let~$f_1(x,y)=f\big(b^{-u}x^q,b^q x^p(1+y)\big)$. If~$f_1$ is nonsingular in~$y$, stop; else, go back to step~\ref{alg:Duval_1} with~$f$ replaced with~$f_1$.
\end{enumerate}
\end{framed}
\caption{Computing parametrisations.}\label{alg:Duval}
\end{Algorithm}

The idea is that we use the Newton polygon to determine the valuation of the roots~$y$ of~$f(x,y)=0$, and then view~$f_0$ as the ``leading terms'', the other terms being thought of as higher-order perturbations. After finitely many iterations, the equation obtained will be nonsingular in~$y$, so its roots can be be found by Newton iteration. We thus obtain explicit parametrisations representing the points of~$\tilde C$ above~$x=0$ such that the field of definition of each point is the extension generated by the coefficients of the corresponding parametrisation. Parametrisations for the points above other values of~$x$ can be of course obtained similarly, by shifting the variable~$x$ appropriately.

\begin{rk}\label{rk:Duval_char0}
The only reason why we have assumed that~$K$ has characteristic 0 was to ensure that~$f(x,y) \in K(x)[y]$ splits completely over~$\overline K \{\{x\}\}$. Theorem~\ref{thm:Duval} and Algorithm~\ref{alg:Duval} actually remain valid in positive characteristic~$\pi$ as long as there is no wild ramification, that is to say that none of the places has ramification index divisible by~$\pi$, which is equivalent to having~$\pi \nmid q$ whenever we consider a segment~$pi+qj=r$ of a Newton polygon in step~\ref{alg:Duval_1}. All the algorithms presented in this section therefore remain valid in positive characteristic as long as~$\tilde C$ is at most tamely ramified as a cover of~$\P^1_x$, which in practice means we typically only exclude really small characteristics such as 2, 3, or 5. Furthermore, by checking whether~$\pi \mid q$ during the execution of algorithm~\ref{alg:Duval}, we can reliably detect when this algorithm is going to fail.
\end{rk}

\subsection{Regular differentials and the genus}

Now that we have computed parametrisations representing singular points and points at infinity, we can find a basis of regular differentials on~$\tilde C$. Indeed, it is well-known~\cite[2.9]{HoloDiffs} that for all~$(i,j) \in \Z$ strictly in the interior of the full (as opposed to lower) convex hull of the support of~$f(x,y) = \sum_{i,j} a_{i,j} y^i x^j$, the differential~$\omega_{i,j} = \frac{x^{j-1} y^{i-1}}{\partial f / \partial y} \operatorname{d}\!x$ is regular everywhere except possibly at singular points, and that every regular differential on~$\tilde C$ is a~$K$-linear combination of those. We thus obtain a basis of regular differentials by finding the linear combinations whose expansion along the parametrisations corresponding to singular points do not have any polar part, which amounts to linear algebra over~$K$. In particular, we recover the genus of~$\tilde C$ as the size of this basis.

While there exist more direct ways to compute the genus, having an actual basis of regular differentials is very useful in practice. For example, it makes it possible to test whether the curve is hyperelliptic, and to find an explicit change of variables which puts in in Weierstrass form if it is~\cite{HoeijHyperell}. And if the curve is not hyperelliptic, on can instead compute its canonical image, which provides a way of finding simpler models for curves defined by a complicated, highly-singular equation (for example, this is the approach that we followed in~\cite[3]{SL3}).

\subsection{Riemann-Roch spaces and extra functionalities}

With our parametrisations representing singular points and points at infinity, we can also compute the integral closure 
\[ \O = \{ s \in K(C) \, \vert \, \text{the only poles of } s \text{ are above } x=\infty \} \]
of~$K[x]$ in~$K(C)$ in a similar way to the number field case~\cite[2.4]{GTM193}: for each irreducible~$d(x) \in K[x]$ such that~$d(x)^2 \mid \disc_y f(x,y)$, we construct a local basis by starting with the approximation~$(\omega_j = y_1^{j-1})_{1 \leqslant j \leqslant n}$ where~$y_1=a(x) y$ and~$a(x)$ is the leading coefficient of~$f(x,y) \in K(x)[y]$, and refining it as long as we can find scalars~$\lambda_j \in K[x] / \big(d(x)\big)$ such that~$\frac{\sum_j \lambda_j \omega_j}{d(x)}$ has no polar part when evaluated along the parametrisations representing the points above~$d(x)=0$. We then join these local bases into a~$K[x]$-basis of~$\O$ by computing a Hermite normal form over~$K[x]$.

Thanks to this~$K[x]$-basis of~$\O$, we can check whether~$C$ is geometrically irreducible, by finding which elements of~$\O$ are also regular above~$x=\infty$.

We can also compute Riemann-Roch spaces, since it is easy, given a divisor on~$C$, to compute a ``common denominator''~$d(x) \in K[x]$ such that the corresponding Riemann-Roch space is contained in ~$\frac1{d(x)}\O$.

This makes it possible to find conic models for curves of genus~$0$. If~$K$ is a number field, we can then test whether the curve has a rational point by a constructive version of Hasse-Minkowski, in which case another use of Riemann-Roch provides us with an explicit rational parametrisation of the curve~\cite{Hoeij0}. Riemann-Roch spaces also make it possible to turn curves of genus 1 on which a rational point is known into elliptic curves in Weierstrass form.

Finally, now that we are able to compute Riemann-Roch spaces, we can initialise Makdisi's algorithms so as to compute in the Jacobian of~$\tilde C$. 

We have implemented all the functionalities described in this section in~\cite{gp}. Our code, which compares quite decently to~\cite{Magma}, is available for use in a development branch of~\cite{gp}, which also contains the generalisation of~\cite{Hensel} to~$\Q(t)$ described in Section~\ref{sect:HenselQt}.

\section{Examples}\label{sect:examples}

In order to demonstrate the use of the algorithm described in Section~\ref{sect:HenselQt}, we have computed some division polynomials over~$\Q(t)$. The calculations took place on the~\cite{Plafrim} cluster.

\subsection{Warmup}

As a sanity check, we first used our new algorithm in order to recover an equation for the~$3$-torsion of the elliptic surface~$\E$ defined by
\[ y^2=t(1+2t-t^2)(x^2-1)(t^2x^2-1) \]
that was the object of our attention in~\cite{SL3}. Even though using Makdisi's algorithms on elliptic curves is obviously out-of-proportion, we instantaneously obtained the division polynomial
\[ 3x^8 + 4t(t^2+1)(t^2-2t-1)x^6+6t^4(t^2-2t-1)x^4-t^8(t^2-2t-1)^2 \in \Q(t)[x], \]
which is incomparably simpler than what we obtained in~\cite{SL3} with~\cite{gp}'s \verb?elldivpol? function, and even prettier than the nicest model that we were able to achieve in~\cite{SL3}. To boost, this polynomial reminisces about~$t=0$ and~$t^2-2t-1=0$ being places of bad reduction of~$\E$.

\subsection{A hyperelliptic family}

Encouraged by this first example, we then computed an~$\ell$-division polynomial for~$\ell=3$ of the curve over~$\Q(t)$ of genus~$g=2$ corresponding to the hyperelliptic surface~$H$ defined by the equation
\[ y^2 = x^6 - x^4 + (t-1)(x^2+x). \]
\begin{rk}
The equation~$y^2 = x^6 - x^4 + t(x^2+x)$ would have been more natural, but we shifted the parameter~$t$ so as to have good reduction at~$t=0$. We did the same for the previous example, but the polynomial which we presented there was the un-shifted version.
\end{rk}

We chose to use the auxiliary prime~$p=17$, since having the~$\ell$-torsion defined over~$\Q_{p^a}((t))$ then merely requires~$a=6$; and we computed the~$\ell$-torsion mod~$(p^e,t^h)$ for~$e=48$ and~$h=16$. The computation took 2 minutes, and we obtained an~$\ell$-division polynomial~$R_{H,3}(x,t) \in \Q(t)[x]$ of degree~$\ell^{2g}-1=80$ and whose coefficients have numerators of degree up to 12 and coefficients of up to 27 decimal digits, and common denominator~$d_H(t) = 3^3(t+1)^2$.

This denominator can probably be explained by the fact that~$H$ has bad reduction at~$t=-1$; even though it can be observed that~$d_H(t)$ is not divisible by~$t-1$ whereas~$H$ clearly has bad reduction at~$t=1$ as well.

\subsection{A plane quartic family}\label{sect:divpol_Q}

As a final example, we computed an~$\ell$-division polynomial for~$\ell=2$ of the family~$Q$ of plane quartics of generic genus~$g=3$ defined by the equation
\[ x^4+(2-t)y^4+2x^3+x(x+y)+(t-1)(y+x^2+x) = 0. \]
This time, we took~$p=5$ as it allows~$a=7$, and the accuracy parameters were~$e=h=128$. After one hour and a half, we obtained a division polynomial~$R_{Q,2}(x,t) \in \Q(t)[x]$ of degree~$\ell^{2g}-1=63$ with common denominator~$d_Q(x) = (t-2)(2t-3)^4 d_{22}(t)$ where~$d_{22}(t) \in \Z[t]$ is irreducible of degree 22 and has leading coefficient~$2^{16}$, and whose coefficient numerators have degree up to 54 and coefficients of up to 39 digits.

It should be noted that one of the places of~$\P^1_t$ at which~$Q$ has bad reduction has degree~$14$ over~$\Q$; since this must somehow be reflected in an anomalous behaviour of the specialisation of~$R_{Q,2}$ at this value of~$t$, this explains why the coefficients of~$R_{Q,2}$ are so complicated, and why the~$t$-adic accuracy ($h=128$) required to identify them was so much larger than in the previous example. This in turn explains why this computation took so much longer than the previous one. 

This time, most of the ``geometric content'' of the denominator, that is to say the factors~$(2t-3)^4$ and~$d_{22}(t)$, do not correspond to places of bad reduction of~$Q$ (but~$t-2$ does), and should instead probably be interpreted as values of~$t$ for which the evaluation map~$\alpha \in \Q(t)(\JJ)$ fails to be defined on all the~$2$-torsion points (see Section~\ref{sect:HenselQt} for the definition and context around~$\alpha$). However, it is still interesting to note that in all three examples, the ``arithmetic content'', that is to say the leading coefficient of the common denominator, is a power of~$\ell$.

\begin{rk}
Our calculations rely on~\cite{gp}'s polynomial arithmetic, which unfortunately does not benefit from fast algorithms for multiplication of polynomials of high degree. In view of the high~$t$-adic accuracy that it required, it is likely that the computation of~$R_{Q,2}$ would have been faster if fast polynomial arithmetic had been available.
\end{rk}

\begin{rk}
As explained in the Introduction, our identification of the coefficients of our division polynomials as elements of~$\Q(t)$ from approximations in~$\Q_p[[t]]$ is not rigorous. However, it is easy to convince oneself that these division polynomials are correct beyond reasonable doubt, for example by checking that their at nonzero values of~$t$ of good reduction has Galois group contained in~$\GSp(2g,\ell)$, and that their ramification agrees what is predicted by N\'eron-Ogg-Shafarevich~\cite{NOS}. The geometric interpretation of the ramification of the specialisations of these division polynomials at bad values of~$t$ which we will establish in the next section is  also evidence that their coefficients have been correctly identified.
\end{rk}

\section{Degeneration of Galois representations and their ramification}\label{sect:degen}

Disappointingly, the division polynomials~$R_{H,3}(x,t)$ and~$R_{Q,2}(x,t)$ which we have obtained in the previous Section are so complicated that neither~\cite{Magma} nor our plane curves package presented in Section~\ref{sect:algcurves} are able to determine their genus, let alone compute Riemann-Roch spaces required to use Makdisi's algorithms to work in their Jacobian. As a result, we are unfortunately unable to conclude our calculation of the Galois representations occurring in the \'etale cohomology of the corresponding surfaces.

However, these division polynomials are still very valuable data, in that each of them encodes a family of Galois representations parametrised by~$\P^1_\Q$. Furthermore, these representations are far from trivial, in that they have maximal image. Indeed, one easily checks with~\cite{Magma} that the specialisation of~$R_{H,3}(x,t)$ at a rational value of~$t$ of good reduction of~$H$ (for example, at~$t=0$) has Galois group~$\GSp(4,3)$ over~$\Q$, which proves that~$R_{H,3}(x,t)$ has Galois group~$\GSp(4,3)$ over~$\Q(t)$; therefore, most specialisations of~$R_{H,3}(x,t)$ will have Galois group~$\GSp(4,3)$ by Hilbert irreducibility, so that~$R_{H,3}(x,t)$ may be viewed as a family (in~$t$) of polynomials (in~$x$) with generic Galois group~$\GSp(4,3)$. One similarly checks that~$R_{Q,2}(x,t)$ defines a family of polynomials with generic Galois group~$\GSp(6,2)=\Sp(6,2)$, which happens to be a simple group.

\subsection{Decomposition of the bad places}

It is especially interesting to study how these families of Galois representations degenerate at values of~$t$ which are places of bad reduction of the corresponding curves over~$\Q(t)$.

The bad places of our hyperelliptic family~$H$ defined by
\[ y^2 = x^6 - x^4 + (t-1)(x^2+x) \] are easily determined by examining the discriminant of the right-hand side, and turn out to be~$t=1$,~$t=-1$,~$t=283/256$, and~$t=\infty$.

In order to analyse the degeneration of~$R_{H,3}(x,t)$ at these places, one must not simply substitute these values for~$t$, as this would be as incorrect as trying to understand the decomposition of a prime~$p$ in a number field by factoring a polynomial mod~$p$ without taking into consideration the index of the order attached to this polynomial.
Instead, we must study the factorisation over~$\Q((t))$ of versions of~$R_{H,3}(x,t)$ shifted in such a way that the bad place under consideration in now~$t=0$. In view of~\eqref{eqn:branch_Puiseux}, this is equivalent to determining the ramification in~$t$ and the field of definitions of the points above~$t=0$ of the desingularisation of the curve~$R_{H,3}(x,t)=0$, which we can achieve thanks to our implementation of Duval's method described in Section~\ref{sect:Duval}. We thus obtain the following data:

\begin{Table}
\[ \begin{array}{c|l|l|l}
t & \text{Place decomposition} & \text{Galois group} & \text{Ramification} \\
\hline
1 & \Q(\sqrt3)^1 \cdot \Q(\sqrt{-1})^3 \cdot \big( \Q(\zeta_9)^+(\sqrt{-1})\big)^9 \cdot \big( \Q(\zeta_{36})^+ \big)^3 & (\Z/36\Z)^\times & 2,3 \\
-1 & \Q(\sqrt{-21})^1 \cdot K_6^1 \cdot K_{18}^1 \cdot {K'_{18}}^3 & C_2 \times C_3 \cdot S_3^2 & 2,3,7,11 \\
\frac{283}{256} & \Q(\sqrt{-14})^1 \cdot {K''_{18}}^3 \cdot K_{24}^1 & (C_2 \times C_3 \rtimes S_3) \cdot S_4 & 2,3,7,11 \\
\infty & \Q^2 \cdot \Q^6 \cdot \Q(\sqrt{3})^4 \cdot \Q(\sqrt[4]{12})^4 \cdot \Q(\sqrt[4]{12})^{12} & D_4 & 2,3
\end{array} \] 
\caption{Decomposition of the bad places of~$H$.}
\label{tab:H}
\end{Table}

In this table, the second column shows the decomposition of the place of~$\Q(t)$ in the function field~$\Q(t)[x]/\big(R_{H,3}(x,t)\big)$; for example, there are five places above~$t=\infty$, two with residue field~$\Q$ and respective ramification indices 2 and 6, one with residue field~$\Q(\sqrt3)$ and ramification index~$4$, and two with residue field~$\Q(\sqrt[4]{12})$ and respective ramification indices 4 and 12. The third column shows the Galois group of the compositum of the Galois closures of the residue fields, and the last column lists the prime numbers which ramify in this Galois closure, or, equivalently, in at least one of the residue fields.
Still in this table,~$\Q(\zeta_m)^+$ denotes the intersection of the cyclotomic field~$\Q(\zeta_m)$ with~$\R$, and~$K_d$,~$K'_d$,~$K''_d$, and so on stand for pairwise non-isomorphic number fields of degree~$d$. As for Galois groups,~$C_n$,~$D_{2n}$, and~$S_n$ respectively denote cyclic, dihedral, and symmetric groups, and~$A \cdot B$ stands for a nonsplit group extension with normal subgroup~$A$ and quotient~$B$. For~$t=1$, we have exceptionally expressed the Galois group as~$(\Z/36 \Z)^\times$ instead of~$C_6 \times C_2$ because the Galois closure is the~$36^\text{th}$ cyclotomic field.

We will elucidate the nature of some of these residue fields in Section~\ref{sect:fibres}, where we will also explain the occurrence of each of the ramified primes.

\bigskip

As for our family of quartics~$Q$, the places of bad reduction are~$t=1$,~$t=2$,~$t=\infty$, as well as the place of degree 14 mentioned in the previous Section. The high degree of this last place makes explicit computations with it impractical, so we ignore it from now on. We obtain the following data:

\begin{Table}
\[ \begin{array}{c|l|l|l}
t & \text{Place decomposition} & \text{Galois group} & \text{Ramification} \\
\hline
1 & \Q^1 \cdot \Q^1 \cdot \Q^1 \cdot K_8^1 \cdot K_8^1 \cdot {K'_8}^2 \cdot {K''_8}^2 \cdot K_{12}^1 & C_2^3 \rtimes S_4 & 2,229 \\
2 & \Q^1 \cdot \Q^2 \cdot \Q^4 \cdot \Q^8 \cdot \Q^8 \cdot \Q(\sqrt2)^4 \cdot \Q(\sqrt2,\sqrt{15})^8 & C_2^2 & 2,3,5 \\
\infty & \Q^1 \cdot \Q^2 \cdot \Q^4 \cdot K_3^2 \cdot K_3^4 \cdot K_6^1 \cdot {K'''_{8}}^4 & S_4 \times C_2 & 2,23 \\
\end{array} \] 
\caption{Decomposition of some of the bad places of~$Q$.}
\label{tab:Q}
\end{Table}

\subsection{Visualising ramification on the special fibre}\label{sect:fibres}

We would now like to find a geometric explanation for the ramified primes observed in the previous tables. We will also explain the occurrence of some of the residue fields.

At a place of~$\P^1_t$ of good reduction, so that the fibre of the surface is a nice curve~$F$, the N\'eron-Ogg-Shafarevich criterion~\cite{NOS} would lead us to expect ramification at~$p=\ell$ as well as at the primes of bad reduction of~$F$. By analogy, at a bad place, we would expect ramification at~$p=\ell$ and at the primes~$p$ such that the bad fibre becomes ``even worse''.

More specifically, this bad fibre should be understood as the fibre of a minimal regular model of the surface over~$\Q$, and saying that the fibre becoming even worse mod~$p$ means that the reduction mod~$p$ of this special fibre does not agree with the special fibre of the minimal regular model of the reduction mod~$p$ of the surface. In more colourful language, this could be summarised by saying that along with~$p=\ell$, these are the primes~$p$ such that taking the special fibre of the minimal regular model does not commute with reduction mod~$p$.

\begin{rk}\label{rk:no_Neron}
Instead of looking at special fibres of the minimal regular model, it would also make sense to consider the semistable fibres. We content ourselves with this imprecision, because we are in effect looking at families of curves over the base~$\P^1_\Z$ which has dimension 2 (one geometric dimension and one arithmetic one), so that as far as the author is aware, there is no longer a canonical notion of good (meaning N\'eron) model for the Jacobian.
\end{rk}


\subsubsection{The hyperelliptic surface}

Let us begin with the hyperelliptic surface~$H$.

\paragraph{The fibre at~$t=1$}\mbox{}\\

 The surface~$H$ is not regular above~$t=1$, but in characteristic~$\pi \neq 2$, it becomes regular after one blowup, and its special fibre then consists of two rational curves arranged as shown on Figure~\ref{fig:fib_h_1}:

\begin{Figure}[H]
\begin{center}
\begin{tikzpicture}[scale=3]
\draw[thick,variable=\t,domain=-1.2:1.2,samples=50]
plot ({\t*\t-1},{\t*\t*\t-\t});
\draw[thick] (0,-0.6)--(0,0.6);
\end{tikzpicture}
\end{center}
\caption{The special fibre of~$H$ at~$t=1$ when~$\pi \neq 2$.}
\label{fig:fib_h_1}
\end{Figure}

In contrast, in characteristic~$\pi=2$, it takes many more blowups to obtain a regular model of~$H$ above~$t=1$. This explains the ramification at~$p=2$ observed in Table~\ref{tab:H} for~$t=1$. As for ramification at~$p=3$, it is simply explained by the fact that we are looking at~$3$-torsion.

\paragraph{The fibre at~$t=-1$}\mbox{}\\

For~$t=-1$, in characteristic~$\pi \not \in \{2,7,11\}$, we again obtain a regular surface after one blowup. Its special fibre is made up of an elliptic curve and a rational curve, as shown on Figure~\ref{fig:fib_h_m1}. Our plane curve package described in Section~\ref{sect:algcurves} informs us that over~$\Q$, the elliptic component is the curve of~\cite{LMFDB} label \href{https://www.lmfdb.org/EllipticCurve/Q/176/a/2}{176.a2}, whose conductor~$176=2^4 \cdot 11$. 

\begin{Figure}[H]
\begin{center}
\begin{tikzpicture}[scale=3]
\draw[thick,variable=\t,domain=-0.95:1.3,samples=50]
plot ({\t/2},{sqrt(\t^3+1)/3})
plot ({\t/2},{-sqrt(\t^3+1)/3});
\draw[thick,variable=\t,domain=-0.5:0.5,samples=50]
plot ({\t*\t/6-1/2},{\t/3-\t^3/18});
\draw[thick] (0,-0.6)--(0,0.6);
\end{tikzpicture}
\end{center}
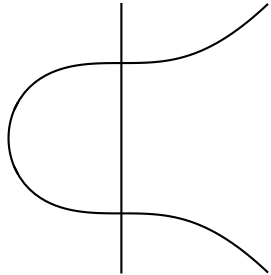
\caption{The special fibre of~$H$ at~$t=-1$ when~$\pi \not \in \{2,7,11\}$.}
\label{fig:fib_h_m1}
\end{Figure}

As a result, in characteristic~$\pi=11$, the elliptic curve degenerates, and the special fibre becomes what is shown on Figure~\ref{fig:fib_h_m1_11}:

\begin{Figure}[H]
\begin{center}
\begin{tikzpicture}[scale=3]
\draw[thick,variable=\t,domain=-1.2:1.2,samples=50]
plot ({\t*\t-1},{\t*\t*\t-\t});
\draw[thick] (-0.5,-0.6)--(-0.5,0.6);
\end{tikzpicture}
\end{center}
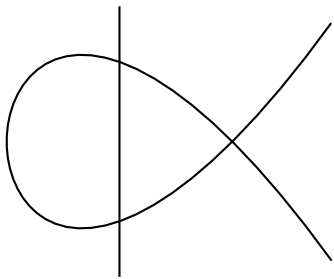
\caption{The special fibre of~$H \bmod 11$ at~$t=-1$. Both components are now rational.}
\label{fig:fib_h_m1_11}
\end{Figure}

This explains why we observed ramification at~$p=11$. As for~$\pi=2$, the special fibre is the same as for~$t=1$, since~$t$ is defined over~$\Z$ and~$-1 \equiv 1 \bmod 2$.

It remains to explain ramification at~$p=7$. A closer inspection of the special fibre over~$\Q$ (as shown on Figure~\ref{fig:fib_h_m1}) shows that the intersection points of the two components are not rational, but defined over~$\Q(\sqrt7)$ and Galois-conjugates of each other; as a result, when we reduce mod~$\pi=7$, these intersection points coalesce, and the special fibre becomes what is shown on Figure~\ref{fig:fib_h_m1_7}, which explains ramification at~$7$:

\begin{Figure}[H]
\begin{center}
\begin{tikzpicture}[scale=3]
\draw[thick,variable=\t,domain=-0.95:1.3,samples=50]
plot ({\t/2},{sqrt(\t^3+1)/3})
plot ({\t/2},{-sqrt(\t^3+1)/3});
\draw[thick,variable=\t,domain=-0.5:0.5,samples=50]
plot ({\t*\t/6-1/2},{\t/3-\t^3/18});
\draw[thick] (-0.5,-0.6)--(-0.5,0.6);
\end{tikzpicture}
\end{center}
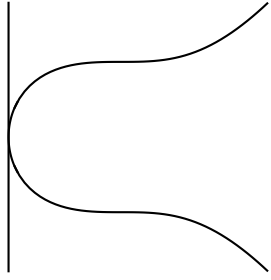
\caption{The special fibre of~$H \bmod 7$ at~$t=-1$.}
\label{fig:fib_h_m1_7}
\end{Figure}

\begin{rk}\label{rk:K18}
As one would expect, our residue fields pick up the~$3$-torsion of the elliptic curve component of the special fibre. More specifically, this elliptic curve \href{https://www.lmfdb.org/EllipticCurve/Q/176/a/2}{176.a2} acquires two of its~$3$-torsion points over~$\Q(\sqrt{-1})$, whereas each of its remaining six points of order~$3$ is defined over one of the Galois conjugates of a number field~$F$ of degree 6. The field~$K_6$ appearing in Table~\ref{tab:H} is actually an extension of~$\Q(\sqrt{-1})$ of degree~$3$ and relative discriminant~$(1+\sqrt{-1})^2 \cdot 3^3 \cdot 7$, whereas the field~$K_{18}$ appearing in the same table is an extension of~$F$ of degree~$3$ ramified only above~$2$ and~$7$. The fact that these extensions have degree 3 can be interpreted in terms of generalised Jacobians, since we are looking at~$3$-torsion. Curiously, there does not seem to be a similar interpretation for~$K'_{18}$, but we still note that~$K_{18}$ and~$K'_{18}$ have the same Galois closure, which also contains~$K_6$.
\end{rk}

\paragraph{The fibre at~$t=\infty$}\mbox{}\\

The surface~$H$ is actually already regular at~$t = \infty$ in any characteristic, so we can directly visualise its special fibre, which turns out to have a rather nasty singularity:

\begin{Figure}[H]
\begin{center}
\begin{tikzpicture}[scale=3]
\draw[thick,variable=\t,domain=-1:1,samples=50]
plot ({\t*\t-1},{\t-2*\t^3+\t^5});
\draw[thick,variable=\t,domain=1:1.3,samples=50]
plot ({\t*\t-1},{\t-2*\t^3+\t^5});
\draw[thick,variable=\t,domain=-1.3:-1,samples=50]
plot ({\t*\t-1},{\t-2*\t^3+\t^5});
\end{tikzpicture}
\end{center}
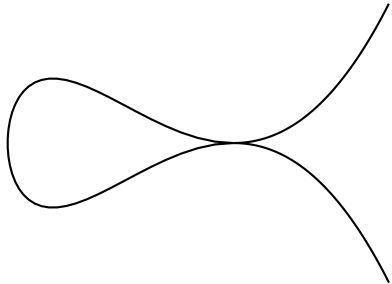
\caption{The special fibre of~$H$ at~$t=\infty$ in any characteristic.}
\label{fig:fib_h_oo}
\end{Figure}

The fact that~$H$ is regular at~$t=\infty$ even mod~$\pi=2$ fails to explain why we observed ramification at~$p=2$ in Table~\ref{tab:H}. However, the special fibre which we have obtained is clearly not semistable, so we may be looking at the ``wrong'' fibre.

In order to investigate further, we can look in the direction of the semistable fibre, which means we must perform a ramified base change~\cite[3.47]{ModuliCurves}. The simplest candidate is to base-change to~$\Q(t^{1/2})$, meaning that we replace~$t$ with~$t^2$ in our equation. This results in~$H$ no longer being regular, even in characteristic~$\pi=0$; after several blowups, we find that in characteristic~$\pi \neq 2$, the special fibre is made up of four rational curves, one of which has multiplicity two, as shown on Figure~\ref{fig:fib_h2_oo}:

\begin{Figure}[H]
\begin{center}
\begin{tikzpicture}[scale=3]
\draw[thick,double,variable=\t,domain=-1:1,samples=50]
plot({\t},{\t*\t-0.5});
\draw[thick,variable=\t,domain=-1:1,samples=50]
plot({\t},{-\t*\t+0.5});
\draw[thick] (0.3,-0.1)--(0.3,-0.7);
\draw[thick] (-0.3,-0.1)--(-0.3,-0.7);
\end{tikzpicture}
\end{center}
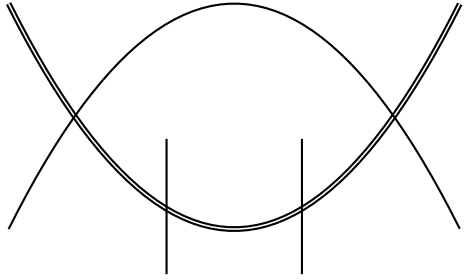
\caption{The special fibre of the base change of~$H$ to~$\Q(t^{1/2})$ at~$t=\infty$ in characteristic~$\pi \neq 2$.}
\label{fig:fib_h2_oo}
\end{Figure}

In contrast, in characteristic~$\pi=2$, the desingularisation requires more blowups, which finally explains the ramification that we observed at~$p=2$.

\begin{rk}
Because of the presence of a double component, the special fibre which we have obtained after base-changing to~$\Q(t^{1/2})$ is still not semistable, and a further base change would be required to remedy this. However, as explained in Remark~\ref{rk:no_Neron}, since we do not have a clear notion of ``good'' model, we content ourself with this reasonably satisfying explanation.
\end{rk}

\paragraph{The fibre at~$t=283/256$}\mbox{}\\

In characteristic~$\pi \not \in \{2,3,7,11\}$,~$H$ is already regular at~$t =283/256$, and its special fibre is a curve of genus 1 with a nodal self-intersection, as shown on Figure~\ref{fig:fib_h_283256}:

\begin{Figure}[H]
\begin{center}
\begin{tikzpicture}[scale=3]
\draw[thick,variable=\t,domain=-1.2:1.2,samples=50]
plot ({\t*\t-1},{\t*\t*\t-\t});
\end{tikzpicture}
\end{center}
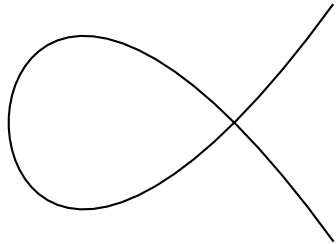
\caption{The special fibre of~$H$ at~$t=283/256$ when~$\pi \not \in \{2,3,7,11\}$.}
\label{fig:fib_h_283256}
\end{Figure}

Over~$\Q$, the desingularisation of this fibre is the elliptic curve of~\cite{LMFDB} label \href{https://www.lmfdb.org/EllipticCurve/Q/528/c/2}{528.c2}, whose conductor is~$528=2^4\cdot 3 \cdot 11$, and as expected, the phenomenon described in Remark~\ref{rk:K18} occurs again, in that the number field~$K_{24}$ displayed in Table~\ref{tab:H} is an extension of degree 3 of the field of degree 8 over the Galois conjugates of which the points of order 3 of this elliptic curve are defined. We do not, however, have a similar interpretation for the field~$K''_{18}$, but we note that its Galois closure is the same as that of~$K_{24}$, and also contains the other residue field~$\Q(\sqrt{-14})$ appearing in the corresponding row of Table~\ref{tab:H}.

In characteristics~$\pi = 2,3,7,11$, we respectively have~$283/256 \equiv \infty,1,-1,-1$, which are cases for which we have already found an explanation for the corresponding ramification.

%
%

%
%

\subsubsection{The quartic surface}

We now proceed to the same analysis of ramification for the family of plane quartics~$Q$.

\paragraph{The fibre at~$t=1$}\mbox{}\\

At~$t=1$, in characteristic~$\pi \not \in \{2,229\}$, we find that the special fibre has three components, two of which are rational, whereas the third one has genus 2:

\begin{Figure}[H]
\begin{center}
\begin{tikzpicture}[scale=3]
\draw[thick,variable=\t,domain=-0.92:1.1,samples=50]
plot ({\t},{sqrt(2*\t^5-\t+1)/3})
plot ({\t},{-sqrt(2*\t^5-\t+1)/3});
\draw[thick,variable=\t,domain=-0.9:0.9,samples=50]
plot ({(\t/3)^2-1},{(\t-10*\t^3/81+40/6561*\t^5)/3});

\draw[thick] (-1,0.5) -- (0.2,-0.2);
\draw[thick] (-1,-0.5) -- (0.2,0.2);
\end{tikzpicture}
\end{center}
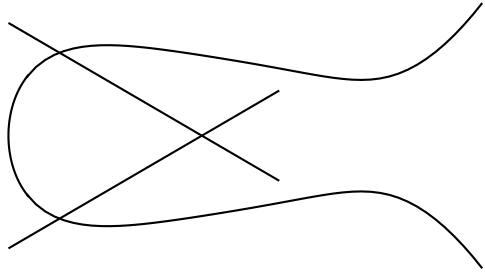
\caption{The special fibre of~$Q$ at~$t=1$ when~$\pi \not \in \{2,229\}$.}
\label{fig:fib_q_1}
\end{Figure}

Over~$\Q$, our plane curves package informs us that the component of genus~2 is isomorphic to the hyperelliptic curve of equation
\[ y^2=x(x^4-x+1) \]
whose~\cite{LMFDB} label is \href{https://www.lmfdb.org/Genus2Curve/Q/29312/a/58624/1}{29312.a.58624.1}; in particular, the conductor of its Jacobian is~$29312=2^7 \cdot 229$. As expected, the phenomenon described in Remark~\ref{rk:K18} occurs again, in that the number field~$K_8$ displayed in Table~\ref{tab:Q} is defined by the irreducible polynomial~$x^8-x^2+1$ and is therefore clearly a quadratic extension of a field over which the Jacobian of this hyperelliptic curve acquires a point of order 2. We do not have any similar interpretation for the fields~$K'_8$,~$K''_8$, nor~$K_{12}$ appearing in the same row of this table, but we still mention that~$K_8$,~$K'_8$, and~$K''_8$ share the same Galois closure, which is a quadratic extension of the Galois closure of~$K_{12}$.

Since~$229$ divides the discriminant of this hyperelliptic curve, when we reduce mod~$\pi=229$, this curve degenerates into a curve of genus 1 with a nodal self-intersection:

\begin{Figure}[H]
\begin{center}
\begin{tikzpicture}[scale=3]
\draw[thick,variable=\t,domain=-1.2:1.2,samples=50]
plot ({(\t*\t-1)*3/2},{(\t*\t*\t-\t)});
\draw[thick] (-1.5,0.5)--(-0.5,-0.15);
\draw[thick] (-1.5,-0.5)--(-0.5,0.15);
\end{tikzpicture}
\end{center}
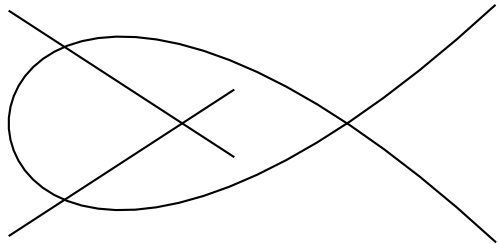
\caption{The special fibre of~$Q \bmod 229$ at~$t=1$.}
\label{fig:fib_q_1_229}
\end{Figure}
This explains the ramification that we have observed at~$p=229$. As for the ramification at~$p=2$, it is explained both by the fact that we are now looking at the~$2$-torsion.

\paragraph{The fibre at~$t=2$}\mbox{}\\

In characteristic~$\pi \not \in \{3,5\}$, we obtain a special fibre made up of three rational components, one of which has a cusp, and which are arranged as follows:

\begin{Figure}[H]
\begin{center}
\begin{tikzpicture}[scale=3]
\draw[thick,variable=\t,domain=-1:1,samples=50]
plot ({\t^3},{\t*\t*2/3});
\draw[thick] (-1,0)--(1,0);
\draw[thick] (0.6,0.7)--(0.6,-0.2);
\end{tikzpicture}
\end{center}
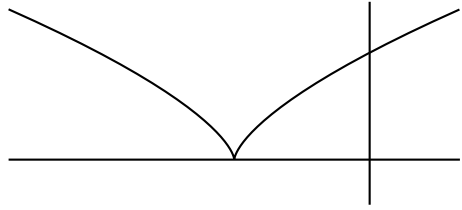
\caption{The special fibre of~$Q$ at~$t=2$ when~$\pi \not \in \{3,5\}$.}
\label{fig:fib_q_2}
\end{Figure}

Reducing mod~$\pi=5$ does not result in requiring more blowups; however, the rightmost fibre, which is a conic, degenerates into a union of two curves, which explains the ramification at~$p=5$: 

\begin{Figure}[H]
\begin{center}
\begin{tikzpicture}[scale=3]
\draw[thick,variable=\t,domain=-1:1,samples=50]
plot ({\t^3},{\t*\t*2/3});
\draw[thick] (-1,0)--(1,0);
\draw[thick] (0.5,0.7)--(1,0.2);
\draw[thick] (0.5,-0.2)--(1,0.35);
\end{tikzpicture}
\end{center}
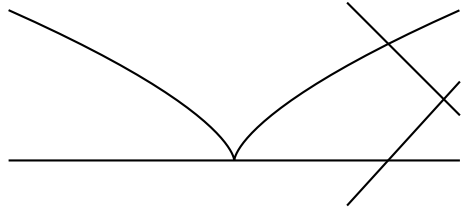
\caption{The special fibre of~$Q \bmod 5$ at~$t=2$.}
\label{fig:fib_q_2_5}
\end{Figure}

The same degeneration occurs mod~$\pi=3$, and furthermore resolving the singularities of~$Q$ at~$t=2$ also requires more blowups in characteristic~$3$. Both these facts explain the ramification at~$p=3$.

\paragraph{The fibre at~$t=\infty$}\mbox{}\\

Mod~$\pi \not \in \{2,23\}$, our model for~$Q$ is already regular at~$t=\infty$, whence a special fibre formed of one component of genus 1 with a nasty self-intersection:

\begin{Figure}[H]
\begin{center}
\begin{tikzpicture}[scale=3]
\draw[thick,variable=\t,domain=-1:1,samples=50]
plot ({\t*\t-1},{\t-2*\t^3+\t^5});
\draw[thick,variable=\t,domain=1:1.3,samples=50]
plot ({\t*\t-1},{\t-2*\t^3+\t^5});
\draw[thick,variable=\t,domain=-1.3:-1,samples=50]
plot ({\t*\t-1},{\t-2*\t^3+\t^5});
\end{tikzpicture}
\end{center}
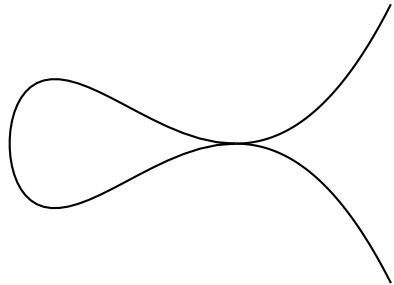
\caption{The special fibre of~$Q$ at~$t=\infty$ when~$\pi \not \in \{2,23\}$.}
\label{fig:fib_q_oo}
\end{Figure}

Over~$\Q$, our plane curve package informs us that the desingularisation of this curve is the elliptic curve with~\cite{LMFDB} label \href{https://www.lmfdb.org/EllipticCurve/Q/92/a/1}{92.a1}, whose conductor is~$92 = 2^2 \cdot 23$; and the field~$K_3$ displayed in Table~\ref{tab:Q}, which is the cubic field of discriminant~$-23$, is also the field over which this elliptic curve acquires a point of order 2. Furthermore,~$K_6$ is a quadratic extension of~$K_3$ which is only ramified above 2 and 23. We do not have a similar explanation for~$K'''_8$, but we observe that the Galois closure of~$K'''_8$, which has degree~$48$, contains~$K_6$ and therefore~$K_3$.

That 23 divides the conductor of this elliptic curve also results in this curve acquiring an extra node mod~$\pi=23$, which explains the ramification at~$p=23$: 

\begin{Figure}[H]
\begin{center}
\begin{tikzpicture}[scale=3]
\draw[thick,variable=\t,domain=-1:1,samples=50]
plot ({\t*\t-1},{-2*\t+5*\t^3-4*\t^5+\t^7});
\draw[thick,variable=\t,domain=-1.5:-1,samples=50]
plot ({\t*\t-1},{-2*\t+5*\t^3-4*\t^5+\t^7});
\draw[thick,variable=\t,domain=1:1.5,samples=50]
plot ({\t*\t-1},{-2*\t+5*\t^3-4*\t^5+\t^7});
\end{tikzpicture}
\end{center}
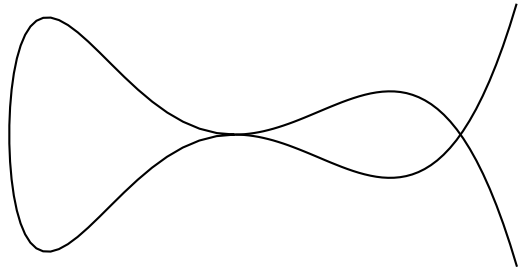
\caption{The special fibre of~$Q \bmod 23$ at~$t=\infty$.}
\label{fig:fib_q_oo_23}
\end{Figure}


\end{document}